\documentclass[12pt]{article}

\usepackage[a4paper, left=3.2cm,top=3.2cm]{geometry}
\usepackage{authblk}
\usepackage[onehalfspacing]{setspace}

\usepackage{amsmath,amssymb,graphicx}
\usepackage{amsmath,amssymb}
\usepackage[colorlinks=true,bookmarks=false,citecolor=blue,urlcolor=blue]{hyperref}
\usepackage[mode=buildnew]{standalone}
\usepackage{tikz}

\newcommand{\R}{\mathbb{R}}
\DeclareMathOperator{\A}{\mathbf{A}}
\DeclareMathOperator{\reg}{\mathbf R}
\DeclareMathOperator{\tik}{\mathbf T}
\DeclareMathOperator{\dist}{\mathbf D}
\newcommand{\norm}[1]{\lVert #1 \rVert}

\allowdisplaybreaks

\title{Relaxed data-consistency for limited bandwidth photoacoustic tomography}

\date{}

\author{Daniel Obmann}

\affil{Department of Mathematics, University of Innsbruck\authorcr
Technikerstrasse 13, 6020 Innsbruck, Austria
 \authorcr E-mail:  \texttt{daniel.obmann@uibk.ac.at}
 }

\author{Markus Haltmeier}

\affil{Department of Mathematics, University of Innsbruck\authorcr
Technikerstrasse 13, 6020 Innsbruck, Austria
 \authorcr E-mail:  \texttt{markus.haltmeier@uibk.ac.at}
 }

\begin{document}

\maketitle

\begin{abstract}
We study the effect of using weaker forms of data-fidelity terms in generalized Tikhonov regularization accounting for model uncertainties. We show that relaxed data-consistency conditions can be beneficial for integrating available prior knowledge.
\end{abstract}

\section{Introduction}
Many optical imaging imaging techniques, including photoacoustic tomography (PAT), can be formulated as an inverse problem of the form $ y  = \A ( x_0 )+ z $, where $\A \in \R^{M  \times N}$ is the forward model describing the specific  imaging system, $y \in \R^M $ are the measured  data, $z$ the noise and $x \in  \R^N$ is the  desired image to be reconstructed.  Common methods for stably computing an approximation $ x_\star \in  \R^N$ of the desired image are based on the following main ingredients: 
\begin{enumerate}
\item \emph{Data consistency:} $\A (x_\star) $ is close to the observed data $y$.      
\item \emph{Prior knowledge:}  $x_\star$ satisfies available  prior information.         
\end{enumerate}   
Tikhonov regularization is an established concept combining data consistency and prior knowledge. It considers minimizers of the Tikhonov functional $\norm{\A (x) - y}^2 + \alpha \reg(x)$ where $\norm{\A (\,\cdot \,) - y}^2$ assures approximate data consistency, $\reg$ is the regularizer that integrates prior knowledge and $\alpha$  is the regularization parameter. 


\begin{figure}[htb!]
    \centering
    \includegraphics[width=0.4\textwidth]{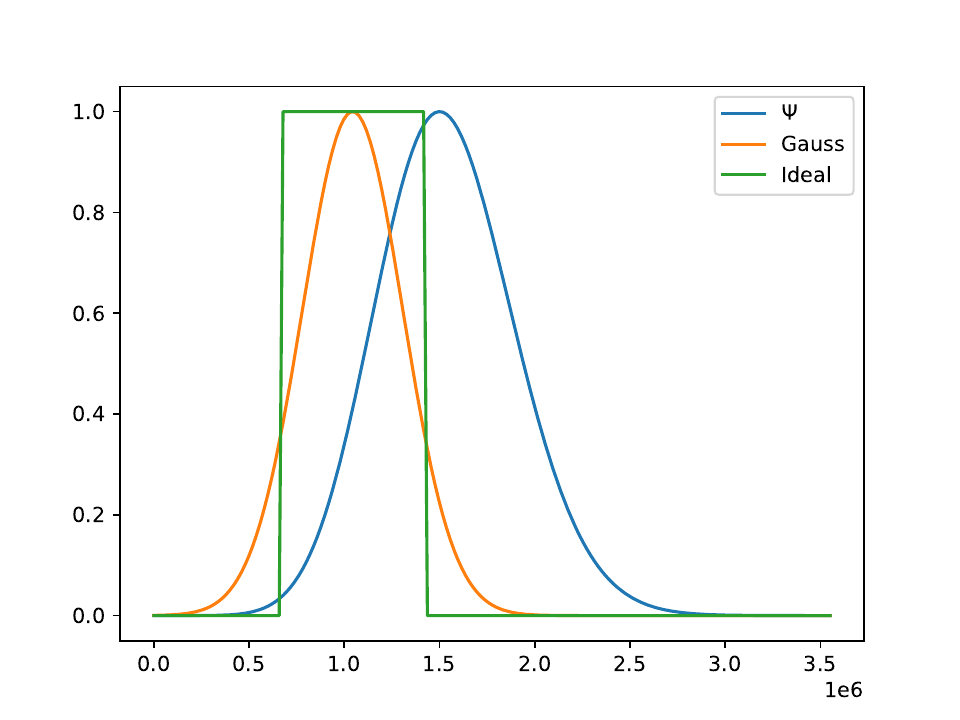}
        \includestandalone{figures/illustration}
    \caption{Left: Frequency response of point spread function $\phi_B$ used to simulate data (blue) and filters $\phi$ used for relaxed data consistency (orange and green). Right: Due to a model mismatch    the prior knowledge  set $C_{\rm prior}$ may not intersect the set of $\epsilon$-data-consistent solutions $\{ x \colon \norm{ \A (x) - y_B } \leq \epsilon \}$ for band-limited data $y_B$. The idea in this paper is to relax the data consistency term such that the true solution $x_0$ becomes consistent with the data.}
    \label{fig:model}
\end{figure}

In practice the imaging model $\A$ may differ from the true underlying physical model $\A_B \in \R^{M  \times N}$ which is either unknown or computationally infeasible.   In this work, we are in particular interested in the case where the available given data $y_B$ are measured using finite bandwidth transducers where the transducers characteristics are not fully known (see Figure \ref{fig:model}, left).  In his case, the band-limited data $y_B$ may actually be far away from $\A( x_0 )$. As a result,  standard data consistency using the approximate forward model $\A$ and the regularizer $\reg$ may become incompatible (Figure \ref{fig:model}, left). That is, whenever data-consistency $\norm{ \A (x) - y_B } \leq \epsilon$ is satisfied, the regularization term is large and vice-versa. Hence minimizers are either not consistent with the given data $y_B$ or they are not regular with respect to $\reg$. 
While one may be tempted to use  a different regularization term, in many  applications we have quite natural priors on the desired image provided  by physical constraints, such as non-negativity of the initial pressure  in PAT.  Band-limiting the data introduces nonphysical negative  components preventing the application of the positivity prior.  In this work we demonstrate that relaxing  data consistency is beneficial  in such scenarios.

\section{Methods}

For the sake of simplicity, we focus on the effect of  transducers bandwidth in 2D PAT. We assume that each transducer has the same frequency response. The measured data can be modeled by $\A_B (x)   =  \phi_B \ast_t \A (x)$, where $\phi_B$ is  the point spread function (PSF) of the detection system  \cite{wang2015photoacoustic}.  Different methods to overcome this problem such as model based \cite{Li:10} and deep learning approaches \cite{awasthi2020deep, gutta2017deep} have been proposed.  For the simulations done here, $\A$ is a discretized PAT forward matrix using the constant sound speed wave equation  using $64$ transducers and $357$ temporal samples.

\subsection{Relaxed data-constancy}

We assume that the exact filter is not available but we have approximate knowledge of the frequencies available in the data.  We  consider minimizing the following  Tikhonov functional with relaxed data-discrepancy term,     
\begin{equation} \label{eq:proposed}
    \tik_{\phi, \alpha} (x)  \triangleq  \dist_\phi( \A (x), y_B  ) + \alpha \reg(x) \,,
    \qquad \text{ with  }  \dist_\phi( \A (x), y_B  ) \triangleq  \norm{  \phi \ast_t ( \A (x) -  y_B ) }^2  \,.
\end{equation} 
Here $\phi$ is a predefined filter, $\norm{\,\cdot\,}$ the standard $\ell^2$-norm, and $\ast_t$ denotes the convolution in the temporal variable.

For the numerical results we test two different instances of the filter  $\phi$. The first one is a Gaussian window centered at some center frequency and the second one is an ideal band-pass filter. We model the filters $\phi$ according to the assumption that we have approximate knowledge and design them to cover these frequencies. An illustration of these filters and the assumed system PSF are shown in Figure~\ref{fig:model} (left).

\subsection{Numerical setup}

To test the effects of relaxing data-consistency we fix the regularizer $\reg$. Here, we choose as a regularization term the total variation (TV) regularization. Additionally, in PAT we have the natural constraint of positivity of the initial pressure. This leads to the regularizer $\reg(x) = \norm{x}_\text{TV} + I_{\geq 0}(x)$ where $\norm{x}_\text{TV}$ it the total variation semi-norm and $I_{\geq 0}$ is the indicator function of the set of non-negative images. We add  Gaussian white noise to  $\A_B (x_0)$ with standard deviation equal to two times the mean of $\A_B (x_0)$. We compare the data fidelity term $\dist_\phi( \A (x), y_B  )$ to the standard squared $\ell^2$-norm distance  $ \dist_0( \A (x), y_B  ) = \norm{  \A (x) -  y_B  }^2$. The regularization parameter is the same for all methods tested and is chosen as $\alpha = 6 \cdot 10^{-6}$. Minimization of each functional is performed by applying $5000$ steps of Chambolle-Pock algorithm \cite{chambolle2011first} as discussed in \cite{sidky2012convex} (Algorithm 4). More details about the measurement geometry, definition of $\A_B$, implementation and code for the presented results are available at \href{https://git.uibk.ac.at/c7021101/pat-proceeding.git}{https://git.uibk.ac.at/c7021101/pat-proceeding.git}.

\begin{figure}[htb!]
    \centering
    \includegraphics[scale=0.8]{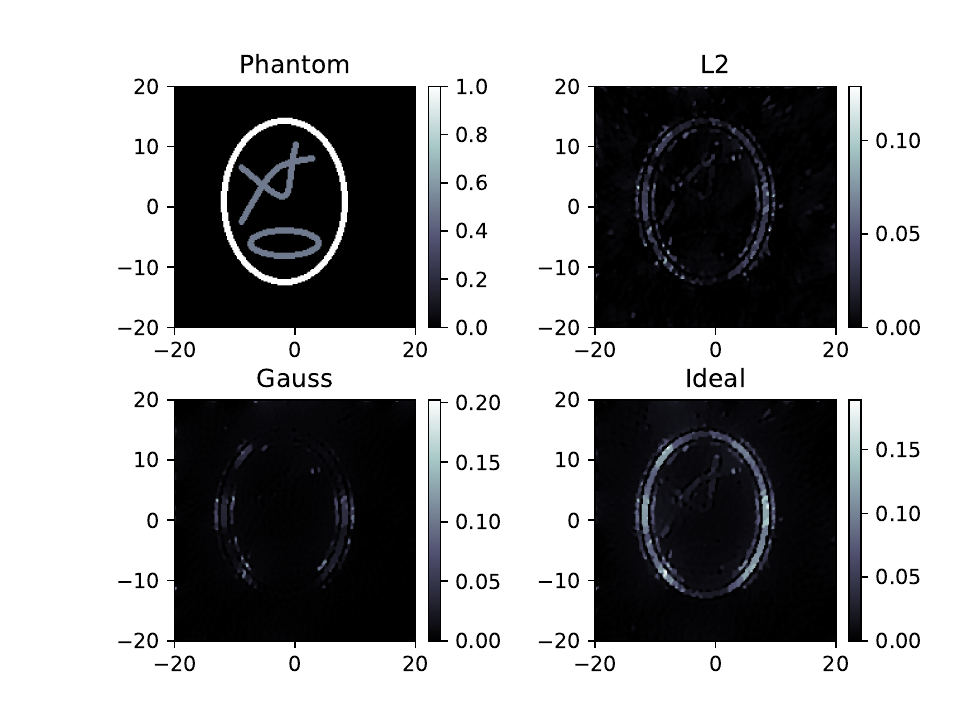}
    \caption{Phantom (top left) used for simulations, $\ell^2$-reconstruction (top right) and reconstructions using the proposed method with Gauss (bottom left) and ideal (bottom right) filter. The reconstructions using the proposed method with an ideal band-pass filter yield higher contrast images and structures are more clearly visible.}
    \label{fig:reconstructions}
\end{figure}

\section{Results}

Reconstruction results with different data-consistency terms are shown in Figure~\ref{fig:reconstructions}. A visual comparison of the three reconstructions shows clear differences. First the reconstruction using $\dist_0$ is able to reconstruct the boundary ellipse and the upside down letter A. However, the boundary ellipse is seemingly disconnected in certain areas and the upside down letter A also shows disconnections and at times even missing structures.
The reconstructions using $\dist_\phi$ with the Gaussian filter barely shows any structure and only the boundary ellipse is visible, but completely disconnected at the bottom. While lowering the regularization parameter in this case improved the visibility of some structures it also introduces highly oscillating and undesirable artefacts.
Using the ideal band-pass filter on the other hand greatly improves the reconstruction quality. First, the boundary ellipse is more clearly visible and is not disconnected as in the case with $\dist_0$ and $\dist_\phi$ using the Gaussian window. Second, the complete structure in the upper part is visible and also connected. Lastly, the reconstruction has a better contrast than the other two reconstructions which further increases the visibility of structures in the image. 
It should be noted, that for the ellipse in the bottom part the reconstruction using $\dist_\phi$ with the ideal band-pass filter is comparable to the one using $\dist_0$ as in both cases the ellipse can barely be seen.


\section{Conclusion}

These results indicate that adjusting the  data-consistency term to the specific problem and in particular adapting them to the measurement devices can have a  beneficial impact on the PAT  image quality. The filters for the simulations performed here were chosen manually and show a huge difference in the reconstruction quality. This further indicates that choosing the correct filter $\phi$ is crucial for obtaining the best possible reconstructions. We believe that an automated filter selection or a filter learned from data could be employed to achieve potentially better results and to reduce the error-prone and tedious task of modelling the filters by hand. In the full proceedings we will additionally present results using different regularizers.

\section*{Acknowledgements}
D.O. and M.H. acknowledge support of the Austrian Science Fund (FWF), project P 30747-N32.

\end{document}